\documentclass[12 pt]{amsart}
\usepackage{amssymb}

\usepackage[mathscr]{eucal}
\normalsize
\theoremstyle{plain}
\newtheorem{Thm}{Theorem}

\theoremstyle{plain}

\theoremstyle{definition}

\renewcommand{\Re}{\operatorname{Re}}
\renewcommand{\Im}{\operatorname{Im}}
\renewcommand{\bar}{\overline}
\newcommand{\st}{\,:\,}
\renewcommand{\tilde}{\widetilde}

\newcommand{\C}{\mathbb C}

\newcommand{\R}{\mathbb R}

\newcommand{\dee}{\partial}

\newcommand{\levi}{\mathscr L}

\DeclareMathOperator{\length}{length}

\DeclareMathOperator{\vol}{vol}
\newcommand{\lam}{\lambda}

\newcommand{\eps}{\epsilon}
\newcommand{\eqdef}{\overset{ \text{def} }{=}}

\newcommand{\bndry}{\text{b}}
\newcommand{\w}{\wedge}

\newcommand{\inv}{^{-1}}

\numberwithin{equation}{section}
\begin{document}
\title[Fefferman's hypersurface measure]{A floating body approach to\\
Fefferman's hypersurface measure\\}
\author{David E. Barrett}
\thanks{This material is based upon work supported in part by the National Science
Foundation under Grant No.
DMS-0072237.}
\address{Dept. of Mathematics\\University of Michigan
\\Ann Arbor, MI  48109-1109  USA }
\date{\today}
\email{barrett@umich.edu }
\begin{abstract} The floating body approach to affine surface area is adapted
to a holomorphic context providing an alternate approach to Fefferman's
invariant hypersurface measure.
\end{abstract}

\maketitle

\section{Introduction}

In [Fef, p. 259] Fefferman introduced a measure $\sigma_Z$ on an arbitrary
smooth strictly pseudoconvex hypersurface $Z$ in $\C^n$.  Viewing
$\sigma_Z$ as a positive
$(2n-1)$-form,  it is characterized  by the
equation 
\begin{equation}\label{E:fefdef}
\sigma_Z\w d\rho =  2^{ 2n/(n+1)} \,M(\rho)^{1/(n+1)}
\omega_{\C^n}
\end{equation}
where $\omega_{\C^n}$ is the euclidean volume form, $\rho$ is a defining
function for $Z$ (i.e., $Z$ is the zero set of $\rho$ and the derivative
of $\rho$ is positive on vectors transverse to $Z$ and pointing to the
pseudoconcave side of $Z$), and
$M$ denotes the complex Monge-Amp\`ere operator defined by
\begin{equation}\label{E:madef}
M(\rho)=(-1)^{n} \det
\begin{pmatrix}
\rho & \rho_{z_j}\\
\rho_{z_{\bar k}} & \rho_{z_{j}\bar z_{  k}}
\end{pmatrix}.
\end{equation}
(The subscripts denote differentiation.)  

The interest in $\sigma_Z$ stems in part from the transformation law 
\begin{equation}\label{E:ftl}
G^* \sigma_{G(Z)} =\left| \det G' \right|^{2n/(n+1)} \sigma_Z 
\end{equation}
valid for $G$ biholomorphic near $Z$ (or for $G$ a CR diffeomorphism
on $Z$).

In the case of a tube hypersurface $Z=X\times i\R^n\subset \R^n\times
i\R^n=\C^n$ it is easy to check (see \S \ref{S:fm} below) that 
\begin{equation}\label{E:fb}
 \sigma_Z= \kappa_{X}^{1/(n+1)}  s_{X} \cdot\omega_{i\R^n};
\end{equation}
here $\omega_{i\R^n}$ is the euclidean volume form on $i\R^n$, $s_{X}$
is euclidean surface area on $X$, and  $\kappa_{X}$ is the Gaussian
curvature of $X$.

The factor $\tilde\sigma_{X}\eqdef\kappa_{X}^{1/(n+1)}  s_{X}$ above
defines a measure on
$X$ which has a longer history; it is the ``affine surface measure"
studied by  Blaschke [Bla].  It satisfies the transformation law
\begin{equation}\label{E:astl}
F^* \tilde\sigma_{F(X)} =\left| \det F' \right|^{(n-1)/(n+1)} \tilde\sigma_X 
\end{equation}
for $F$ affine.

In the case of $\R^2$ Blaschke provided an alternate description  which
applies to general convex curves.  In recent years several works have provided
similar results in higher dimensions.  (For an overview see
[Lei1].) Some of these approaches do not seem to lend
themselves to natural generalization to several complex variables, but one
approach is promising for this purpose, namely that taken in  papers by
Leichtwei{\ss} \,[Lei2] and by Sch\"utt and Werner [ScWe] 
using ``floating body" theory, building on earlier work of Blaschke.

A convex body $K\subset\subset \R^n$  and a positive quantity $\delta$
determine a  {\em convex floating body}  defined to be the
intersection of all
closed half-spaces
$H$ such that $K\setminus H$ has volume $\delta$.  It is common to denote
this object by $K_\delta$, but for notational convenience in this paper we will
let $K_\delta$ denote the portion of $K$ lying outside the convex floating
body.

For $n=3$ and $K$ strictly convex with analytic boundary, Blaschke showed
[Bla] that the affine surface area of $\bndry K$ coincides with
\begin{equation}\label{E:bt}
\sqrt{\pi} \lim_{\delta\searrow0}\frac{\vol( K_\delta)}{\sqrt{\delta}}.
\end{equation}
For general $n$ and $K$ strictly convex with $C^2$ boundary, Leichtwei{\ss}
showed [Lei2]  that the affine surface area of
$\bndry K$ coincides with
\begin{equation}\label{E:leilim}
\lim_{\delta\searrow0} c_n 
\frac{\vol( K_\delta)}{\delta^{2/(n+1)}},
\end{equation} where
\begin{equation*}
c_n=\frac{(2\pi)^{(n-1)/(n+1)}}
{\left(\Gamma\left(\frac{n+1}{2}\right)\right)^{2/(n+1)}}.
\end{equation*} 

In [ScWe] it is shown that for any bounded convex body $K$ in $\R^n$  the
limit
\begin{equation}\label{E:swlim}
\lim_{\delta\searrow0} c_n 
\frac{\vol( K_\delta)}{\delta^{2/(n+1)}}
\end{equation}
 exists and is finite, coinciding with the affine surface area
whenever
$K$ has
$C^2$ boundary. (See \S\ref{E:com} below for more on this result.)

In this paper we provide a generalization of the results \eqref{E:bt} and
\eqref{E:leilim} to Fefferman's measure.

Let $\Omega\subset \C^n$ be a bounded strictly pseudoconvex domain with
$C^{3}$  boundary.   For $M>0$ let
$P_M(\Omega)$ denote the set of
$C^{3}$ functions $h$ on $\bar\Omega$ satisfying the conditions
\begin{enumerate}
\item $h$ is holomorphic on $\Omega$;
\item $h$ and all its derivatives of order $\le 3$ are bounded in
absolute value by $M$ on $\bar\Omega$;
\item  $\bar\Omega\cap h\inv(0)$ is a non-empty subset of
$\bndry\Omega$;
\item $|dh|\ge M\inv$ on  $\bar\Omega\cap h\inv(0)$.
\end{enumerate}

Note that while $P_M(\Omega)$ is not in general biholomorphically invariant,
if $G:\bar\Omega_1\to\bar\Omega_2$ is a $C^{3}$ diffeomorphism
holomorphic on $\Omega_1$ then for $M>0$ there are
$M_\sharp>M_\flat>0$  so that $P_{M_\flat}(\Omega_2)\circ G \subset
P_M(\Omega_1)
\subset P_{M_\sharp}(\Omega_2)\circ G$.

For $\delta>0$ let 
\begin{equation}\label{E:cfldef}
\Omega_{M,\delta}
=\bigcup\limits_{h\in P_M(\Omega)}
\{ z\in \Omega \st 
\vol
\left(
\{w\in\Omega\st |h(w)|\le|h(z)|\}
 \right)< \delta\}. 
\end{equation}

\begin{Thm}\label{T:mt}
For $\Omega$ as above and for all $M\ge M_0(\Omega)$  we have
\begin{equation}\label{E:mt}
C_n
\lim_{\delta\searrow0}
\frac{\vol( \Omega_{M,\delta})}{\delta^{1/(n+1)}}=\int_{\bndry\Omega}
\sigma_{\bndry\Omega},
\end{equation}
where $C_n$ denotes the constant
\[
\left(
\frac{2^{2n-2}\pi^{n-\frac12}\Gamma(\frac{n}{2})}
{(n+1)\Gamma(\frac{n+1}{2})\Gamma(n)}
\right)^\frac{1}{n+1}.\]
\end{Thm}

This theorem will be proved in \S \ref{S:pmt}.  \S \ref{S:fm} has more
information concerning the construction of Fefferman's measure.    The final
section lays out some open questions.

It may strike some readers at this stage that when generalizing the floating
body construction to the holomorphic setting it would seem natural to focus on
sublevel sets of $\Re h$ rather than $|h|$.  Let us address this first in the
one-dimensional setting, where  $\sigma_{\bndry\Omega}$ is  the standard element of arc
length
$|dz|$ (see \S \ref{S:fm} below).  Attempts to understand the volume of small
sublevel sets
$\Re h$ lead to consideration of second-order information about
$\bndry\Omega$ -- but such information simply doesn't appear in the integral
$\int_{\bndry\Omega}
\sigma_{\bndry\Omega}=\length(\bndry\Omega)$.  But for $h\in
P_M(\Omega)$ and $\eps>0$ small the set $\Omega\cap |h|\inv ([0, \eps])$ is
approximately a half-disk, and the parameter $M$ gives us enough uniformity
to assert that for small $\delta>0$ the set $\Omega_{M,\delta}$ is  a
collar about $\bndry\Omega$ of normal width approximately
$\sqrt{2\delta/\pi}$, hence \[\sqrt{\pi/2}\,\lim_{\delta\searrow 0}
\frac{\vol(
\Omega_{M,\delta})}{\sqrt{\delta}}=\length(\bndry\Omega)\] as
claimed in the theorem. 

In higher dimensions the focus on $|h|$ rather than $\Re h$ allows us to
restrict our consideration of second-order information to the complex
directions in the tangent spaces of $\bndry\Omega$.  A related point is that
the small sublevel sets of $|h|$ reflect the non-isotropic structure of
$\bndry\Omega$ (see for example [Rud, \S 5.1]).

\section{On the construction of Fefferman's measure}\label{S:fm}

Let $Z$ be a $C^2$ strictly pseudoconvex hypersurface in an
$n$-dimensional complex manifold $M$ equipped with a smooth positive $2n$-form
$\omega$.  We will explain how to construct a positive $(2n-1)$-form
$\sigma_{Z,\omega}$ on
$Z$ in such a way that the transformation law
\begin{equation*}
G^* \sigma_{G(Z),\tilde\omega} =\left( \frac{G^*\tilde\omega}{\omega}
\right)^{n/(n+1)}
\sigma_{Z,\omega} 
\end{equation*}
holds for $G$ biholomorphic.

Let $J$ denote the complex structure tensor (thus in $\C^n$ we have
$J\frac{\dee}{\dee x_j}=\frac{\dee}{\dee y_j},$ $J\frac{\dee}{\dee
y_j}=-\frac{\dee}{\dee x_j}$).

The Levi-form $\levi$  of $Z$ may be naturally defined as a
symmetric
$TM/TZ$-valued form on $TZ\cap JTZ$ characterized by the identity
\begin{equation*}
\levi(Y_1,Y_2)\equiv [Y_1,JY_2] \mod TZ\cap JTZ
\end{equation*}
for $TZ\cap JTZ$-valued vector fields $Y_1$ and $Y_2$.  The Levi-form is
(real-)hermitian (i.e., 
$\levi(JY_1,JY_2)=\levi(Y_1,Y_2)\,$).  (The hermitian property follows directly
from the integrability condition
$J\left([Y_1,Y_2]-[JY_1,JY_2]\right)=[JY_1,Y_2]+[Y_1,JY_2]$; the symmetry of
$\levi$ follows from the hermitian property and the antisymmetry of the
bracket operation.) Note also that $\levi(Y,JY)=0$.

We carry out the construction first in the two-dimensional case.  

Let $Y$ be a non-zero vector in $TZ\cap JTZ$.  Then $Y, JY, \levi(Y,Y)$ gives
a basis for $TZ$.  We describe $\sigma_{Z,\omega}$ by the identity
\begin{equation}\label{E:fm2}
\sigma_{Z,\omega}\left(Y,JY,\levi(Y,Y)\right)
=\omega^{2/3}\left(Y,JY,\levi(Y,Y), J\levi(Y,Y)\right).
\end{equation}
(We assume here that orientations have been chosen so that $Y,JY,\levi(Y,Y)$
and
\linebreak
$Y,JY,\levi(Y,Y), J\levi(Y,Y)$ are positive bases for $TZ$ and $TM$
respectively.)

If $Y$ is replaced by $\tilde Y=\alpha Y + \beta JY$ then both sides of
\eqref{E:fm2} pick up a factor of $(\alpha^2+\beta^2)^2$; it follows that
$\sigma_{Z,\omega}$ does not depend on the choice of $Y$.

In higher dimension we  choose a complex basis $Y_1,\dots Y_{n-1}$ of
$TZ\cap JTZ$.  Let $L_{j,k}=\levi(Y_j,Y_k)-i\,\levi(JY_j,Y_k)$.  (Thus
$\Big(L_{j,k}\Big)$ is the (complex-)hermitian matrix representing
$\levi$ with respect to the given basis.)  Using the Levi-form to
orient
$TM/TZ$, note that
${\det}^{1/(n-1)}\Big(L_{j,k}\Big)$ defines a vector in 
$TM/TZ$.
  We then describe
$\sigma_{Z,\omega}$ by the identity
\begin{multline}\label{E:fmh}
\sigma_{Z,\omega}\bigg(Y_1,JY_1,\dots,
Y_{n-1},JY_{n-1},{\det}^{1/(n-1)}\Big(L_{j,k}\Big)\bigg)
\\=\omega^{n/(n+1)}\bigg(Y_1,JY_1,\dots,
Y_{n-1},JY_{n-1},{\det}^{1/(n-1)}\Big(L_{j,k}\Big),
J{\det}^{1/(n-1)}\Big(L_{j,k}\Big)\bigg).
\end{multline}

If $Y_1,\dots, Y_{n-1}$ are replaced by $\sum \alpha_{1,k} Y_k + \sum
\beta_{1,k} JY_k,\dots,\sum \alpha_{n-1,k} Y_k + \sum
\beta_{n-1,k} JY_k$, (with $\alpha_{j,k}$ and $\beta_{j,k}$ real) then both
sides of \eqref{E:fmh} pick up a factor of
$\left|\det(\alpha+i\beta)\right|^{2n/(n-1)}$; as before if follows that
$\sigma_{Z,\omega}$ does not depend on the choice of $Y_1,\dots Y_{n-1}$.

We claim that for $M=\C^n$ equipped with the euclidean volume
form $\omega$ the form $\sigma_{Z,\omega}$ defined in \eqref{E:fmh}
coincides with the form $\sigma_Z$ defined in 
\eqref{E:madef}. It will suffice to check this at the origin under that
assumption that $\rho$ is locally of the form $\psi(z_1,\dots,z_{n-1},\Re
z_n)-\Im
z_n$ with $\psi$ and its gradient vanishing at
$0$.  Then 
\begin{equation}\label{E:fg}
\sigma_Z= 2^{\frac{2n}{n+1}} M(\rho)^{\frac{1}{n+1}}  dx_1\w dy_1\dots \w
dx_{n-1}
\w dy_{n-1} \w dx_n 
\end{equation} 
at $0$; setting $Y_j=4\Re\left(\frac{\dee\rho}{\dee
z_n}\frac{\dee}{\dee z_j}-\frac{\dee\rho}{\dee z_j}\frac{\dee}{\dee
z_n}\right), 1\le j\le n-1,$ in \eqref{E:fmh} and checking that $L_{j,k}=
4\rho_{z_j, \bar z_k}\cdot\frac{\dee}{\dee x_n}$ and  
${\det}\Big(L_{j,k}\Big)=2^{2n}M(\rho)\left(\frac{\dee}{\dee
x_n}\right)^{n-1}$at
$0$ we have
\[
2^{\frac{2n}{n-1}}
M(\rho)^{\frac{1}{n-1}}\,\sigma_{Z,\omega}\left(\frac{\dee}{\dee
x_1},\frac{\dee}{\dee
y_1},\dots,\frac{\dee}{\dee x_{n-1}},\frac{\dee}{\dee
y_{n-1}},\frac{\dee}{\dee
x_n}\right)=2^{\frac{4n}{n-1}\frac{n}{n+1}}
M(\rho)^{\frac{2}{n-1}\frac{n}{n+1}}
\]
at $0$;  with a little further manipulation of
exponents we find that $\sigma_{Z,\omega}=\sigma_Z$ as claimed.

Note that using \eqref{E:fg} we may write $\sigma_Z$ in completely euclidean
terms -- up to a multiplicative constant -- as $|\det \levi|^{1/(n+1)}\, s_Z$,
where
$s_Z$ is the euclidean surface area on $Z$ and the bars on $|\det \levi|$
indicate measurement with respect to the euclidean structure.  In the case of a
tube domain $Z=X\times i\R^n$,
$\levi$ is essentially just the second fundamental form of $X$, so
$|\det
\levi|$ is just the Gaussian curvature of $X$.  To check that no
multiplicative constant is missing from
\eqref{E:fb} one can trace through the construction or test both sides
against the hypersurface $\{(z_1,\dots,z_n)\st x_1^2+\dots x_n^2=1\}$.
 ({\em Remark:}  In [Fef], Fefferman allows a
dimension-dependent constant factor in the definition of
$\sigma_Z$; we have chosen the constant $2^{ (2n+1)/(n+1)}$ in
\eqref{E:fefdef} to arrange that
\eqref{E:fb} holds.   A different choice
appears in [Hir]. )

\bigskip

For $n=1$ many of the above computations are problematic, but  we can see
that in this case the natural analogue of the above construction  is given by the formula
\begin{equation*}
\sigma_{Z,\omega}\left(Y\right)
= \omega^{1/2}\left(Y,JY\right)
\end{equation*}
converting a positive area form on $M$ to a positive one-form on
$Z$.  In particular, for $M=\C$ equipped with the euclidean area
form $\omega$, $\sigma_{Z,\omega}$ is  the standard arc
length form, agreeing with $\sigma_Z$ given by \eqref{E:fefdef}.

\section{Proof of main theorem}\label{S:pmt}

Fix for the moment a function $h\in P_M(\Omega)$ and a point $p\in
\bndry\Omega$ where $h$ vanishes.

Choose a unitary system of coordinates $\left(w_1,\dots, w_n\right)$
vanishing at $p$ so that the tangent space to $\bndry\Omega$ is given by
$\Im w_n=0$.  Since the zero set of $h$ must be tangent to $\bndry\Omega$,
we have $dh=h_{w_n}\,dw_n$ at $0$.  Replacing $h$ by
$h_{w_n}\inv(0)\cdot h$ we may assume that $dh=dw_n$ at $0$, this at the
cost of squaring
$M$.

Let $\gamma$ be the local solution to the ordinary differential equation
\begin{equation*}
\gamma_{w_1} h_{w_n} - \gamma_{w_n} h_{w_1}=1
\end{equation*}
subject to the initial condition
\begin{equation*}
\gamma(0,w_2,\dots,w_n)=0.
\end{equation*}

Then the functions $z_1,\dots,z_n$ defined by 
\begin{align*}
z_1&=\gamma(w)\\
z_j&=w_j \text{ for } 2\le j\le n-1\\
z_n&=h(w)
\end{align*}
define a volume-preserving holomorphic change of coordinates near
$p$. 

Note that $\frac{\dee z_j}{\dee w_k}(0)=\delta_{j,k}$; thus this change
of coordinates preserves distances up to a factor of $1+O(\|z\|)$.

With this set up we wish to study the volumes of the  sets
\begin{equation*}
S_\eta\eqdef \{z\in\Omega\st |z_n|<\eta\}
\end{equation*}
where $\Omega$ is locally described by an inequality
\begin{equation*}
\Im z_n > \psi(z_1,\dots,z_{n-1},\Re z_n) 
\end{equation*}
and satisfies 
\begin{equation}\label{E:nocut}
z_n\ne 0 \text{ in }\Omega.
\end{equation}

Let us focus for the time being on the case $n=2$.  Set $ z_1=z, z_2=u+iv$.
Then we may write
\begin{equation}\label{E:psi2}
\psi(z,u)= \lam(u) |z|^2 + \Re \mu(u)z^2+\Re \nu(u) z + \xi(u)+ O(|z|^3) 
\end{equation}
with  $\nu(0)=\xi(0)=\xi'(0)=0$.

The strict pseudoconvexity of $\Omega$ implies that $\lam(0)>0$ and
condition
\eqref{E:nocut} implies that $|\mu(0)|\le\lam(0)$

Let $\tilde\lam=\sqrt{\lam^2+|\mu|^2}$.  Note that
$2|\mu(0)|^2\le\lam(0)^2+|\mu(0)|^2=\tilde\lam^2(0)$, so 
\begin{equation}\label{E:lmi}
|\mu(0)|\le\frac{1}{\sqrt{2}} \tilde\lam(0)
\end{equation}

Let 
\begin{align}\label{E:tpdef}
\tilde\psi(z,u)&=\psi(z,u)+\left(\tilde\lam(u)-\lam(u)\right)|z|^2\notag\\
&= \tilde \lam(u) |z|^2 + \Re \mu(u)z^2+\Re \nu(u) z + \xi(u)+
O(|z|^3),
\end{align}
and let $\tilde\Omega\subset\Omega$ be a domain defined 
near $p$ by the inequality
$v>\tilde\psi(z,u)$. 

For $\eta$ small we may use \eqref{E:lmi} and \eqref{E:tpdef} to conclude that
on
$\bndry\tilde\Omega\cap S_\eta$ we have
\begin{align*}
\left(1-\frac{1}{\sqrt{2}}\right)\tilde\lam(0)|z|^2
&\le \tilde\lam(0)|z|^2+\Re\mu(0)z^2\\
&=v+O\left(\eta(\eta+|z|)\right)+O(|z|^3);
\end{align*}
thus 
$
|z|^2=O(\eta(1+|z|))
$
and $
|z|=O(\sqrt{\eta}),
$
so $\tilde\lam(0)|z|^2+\Re\mu(0)z^2=v+O(\eta^{3/2}).$

Quoting the
fact that $\{z\st A|z|^2+\Re Bz^2< V\}$ has area equal to
$\dfrac{\pi V^+}{\sqrt{A^2-|B|^2}}$ when $|B|<A$  we find that
$\left\{z
\st (z,u+iv)\in \tilde\Omega\right\}$ has area equal to \[\dfrac{\pi
v^+ +O(\eta^{3/2})}{\sqrt{\tilde\lam^2(0)-|\mu(0)|^2}}=\dfrac{\pi
v^+ +O(\eta^{3/2})}{\lam(0)}.\]

Thus
\begin{align}\label{E:vlb}
\vol\left(\left\{ (z,u+iv)\in\Omega\st
|u+iv|<\eta\right\}\right) 
&\ge
\vol\left(\left\{ (z,u+iv)\in\tilde\Omega\st
|u+iv|<\eta\right\}\right)\notag\\
&=\iint\limits_{u^2+v^2<\eta^2} 
\dfrac{\pi
v^+ +O(\eta^{3/2})}{\lam(0)}\,du\,dv\notag\\
&=\int\limits_0^{2\pi} \int\limits_0^\eta 
\dfrac{\pi r \sin^+ \theta+O(\eta^{3/2})}{\lam(0)}\,r\,dr\,d\theta\\
&= \dfrac{2\pi \eta^3+O(\eta^{7/2})}{3\lam(0)}\notag\\
&= \dfrac{8\pi
\eta^3+O(\eta^{7/2})}{3\left(\frac{\sigma_Z}{s_Z}(0)\right)^3}.
\notag
\end{align}

The above estimates are uniform in $p$ and show that $\Omega_{M,\delta}$
is contained in a collar about $\bndry\Omega$ of normal thickness
$
\left(\left(\dfrac{3\delta}{8\pi}\right)^{1/3}+O\left(\delta^{1/2}\right)\right)
\dfrac{\sigma_Z}{s_Z}(0).
$

Thus 
\begin{equation}\label{E:2ls}
\left(\frac{ 8\pi}{3}\right)^{1/3}\limsup_{\delta\searrow0}
\frac{\vol( \Omega_{M,\delta})}{\delta^{1/3}}\le\int_{\bndry\Omega}
\sigma_{\bndry\Omega}.
\end{equation}

To get an estimate in the other direction we make use of that fact that when
$M$ is large enough, for each $p\in\bndry \Omega$ we can find $h_p\in
P_M(\Omega)$ such that
\begin{itemize}
\item $\bar\Omega\cap h_p\inv(0)=\{p\}$;
\item $\|dh_p(p)\|=1$;
\item $\Omega\cap |h_p|\inv ([0, \eps])=\{z\in\Omega\st |h_p(z)|\le\eps\}$ is
connected for $\eps<\eps_0$ (with $\eps_0$ independent of $p$);
\item after introducing new coordinates as above we have $\mu(0)=0$.
\end{itemize}
(See [HeLe, \S 2.4], [Kra, \S 5.2].)   Then
\eqref{E:vlb} can be revised to read
\begin{equation*}
\vol\left(\left\{ (z,u+iv)\in\Omega\st
|z|<\eta\right\}\right) 
= \dfrac{8\pi
\eta^3+O(\eta^{7/2})}{3\left(\frac{\sigma_Z}{s_Z}(0)\right)^3}. 
\end{equation*}

As above, it follows that $\Omega_{M,\delta}$
contains a collar about $\bndry\Omega$ of normal thickness \linebreak
$
\left(\left(\dfrac{3\delta}{8\pi}\right)^{1/3}+O\left(\delta^{1/2}\right)\right)
\dfrac{\sigma_Z}{s_Z}(0),
$
implying that 
\begin{equation}\label{E:2li}
\left(\frac{ 8\pi}{3}\right)^{1/3}\liminf_{\delta\searrow0}
\frac{\vol( \Omega_{M,\delta})}{\delta^{1/3}}\ge\int_{\bndry\Omega}
\sigma_{\bndry\Omega}.
\end{equation}

Combining \eqref{E:2ls} and \eqref{E:2li} we have \eqref{E:mt} in the case
$n=2$.

To treat the case $n>2$ we modify the argument as follows. We now set
$(z_1,\dots,z_{n-1})=z', z_n=u+iv$.
The expansion \eqref{E:psi2} now reads
\begin{equation*}
\psi(z',u)= \sum_{j,k=1}^{n-1}\lam_{j,k}(u) z_j \bar z_k + \Re
\sum_{j,k=1}^{n-1}\mu_{j,k}(u)z_j z_k+\Re\sum_{j=1}^{n-1}
\nu_j(u) z_j +
\xi(u)+ O(\|z'\|^3) 
\end{equation*}
 with  $\nu_j(0)=\xi(0)=\xi'(0)=0$, $\lam_{k,j}=\bar{\lam_{j,k}}$,
$\mu_{k,j}=\mu_{j,k}$.

We may choose an invertible linear map
$T=(T_1,\dots,T_{n-1}):\C^{n-1}\to\C^{n-1}$  and
$\phi_j\ge0$ so that
\begin{equation*}
\psi(z',0)= \sum_{j=1}^{n-1} \left(\left|T_j z'\right|^2+\Re\phi_j
\left(T_jz'\right)^2\right)+ O(\|z'\|^3). 
\end{equation*}
(See for example Lemma 4.1 in [Web].)

Condition
\eqref{E:nocut} implies that each $\phi_j\le 1$.

Let $\big(\tilde\lam_{j,k}\big)$ be the hermitian matrix satisfying
\begin{equation*}
\sum_{j,k=1}^{n-1}\tilde\lam_{j,k} z_j \bar z_k
=\sum_{j=1}^{n-1}\sqrt{1+\phi_j^2} \left|T_j z'\right|^2.
\end{equation*}
In analogy to \eqref{E:lmi} we have 
\begin{equation}\label{E:lmih}
\Re
\sum_{j,k=1}^{n-1}\mu_{j,k}(0)z_j z_k
\le \frac{1}{\sqrt{2}}
\sum_{j,k=1}^{n-1}\tilde\lam_{j,k} z_j \bar z_k.
\end{equation}

We now let 
\begin{align*}
\tilde\psi(z',u)&=\psi(z',u)+\sum_{j,k=1}^{n-1}
\left(\tilde\lam_{j,k}-\lam_{j,k}(u)+O(|u|)\right)
z_j \bar z_k\\ 
&= \sum_{j,k=1}^{n-1}\tilde\lam_{j,k} z_j \bar z_k + \Re
\sum_{j,k=1}^{n-1}\mu_{j,k}(u)z_j z_k+\Re\sum_{j=1}^{n-1}
\nu_j(u) z_j +
\xi(u)+ O(\|z'\|^3),
\end{align*}
where the $O(|u|)$ term is chosen so that there is a domain $\tilde\Omega\subset\Omega$ 
defined near $p$ by \linebreak
$v>\tilde\psi(z',u)$. 
On $\bndry\tilde\Omega\cap S_\eta$ we have as before $
\|z'\|=O(\sqrt{\eta}),
$
and
\[
\sum_{j,k=1}^{n-1}\tilde\lam_{j,k} z_j \bar z_k + \Re
\sum_{j,k=1}^{n-1}\mu_{j,k}(0)z_j z_k=v+O(\eta^{3/2}).
\]

The set 
\begin{align*}
&\left\{z'\st \sum_{j,k=1}^{n-1}\tilde\lam_{j,k} z_j \bar z_k + \Re
\sum_{j,k=1}^{n-1}\mu_{j,k}(0)z_j z_k<V\right\}\\
=&\left\{z'\st \sum_{j=1}^{n-1}
\sqrt{1+\phi_j^2}\left|T_j z'\right|^2+\Re \phi_j 
\left(T_j z'\right)^2
<V\right\}
\end{align*} 
has volume
\begin{align*}
\left| \det T\right|^{-2} \,\prod_{j=1}^{n-1}
\frac{1}{\sqrt{1+\phi_j^2}+\phi_j}
\frac{1}{\sqrt{1+\phi_j^2}-\phi_j}
\cdot\vol\left(\left\{z'\st \|z'\|
<V\right\}\right)
= \frac{\pi^{n-1}(V^+)^{n-1}}{(n-1)! \left| \det T\right|^{2}};
\end{align*} 
thus 
\begin{align*}
\vol\left(\left\{ (z',u+iv)\in\Omega\st
|u+iv|<\eta\right\}\right) 
&\ge
\vol\left(\left\{ (z',u+iv)\in\tilde\Omega\st
|u+iv|<\eta\right\}\right)\notag\\
&=\iint\limits_{u^2+v^2<\eta^2} 
\dfrac{\pi^{n-1}
(v^+)^{n-1} +O\left(\eta^{n-\frac{1}{2}}\right)}{(n-1)!\left| \det
T\right|^{2}}\,du\,dv\notag\\ &=\int\limits_0^{2\pi} \int\limits_0^\eta 
\dfrac{\pi^{n-1} r^{n-1} \left(\sin^+
\theta\right)^{n-1}+O\left(\eta^{n-\frac{1}{2}}\right)}{(n-1)!\left| \det
T\right|^{2}}\,r\,dr\,d\theta\\ 
&=
\dfrac{\pi^{n-\frac12}\eta^{n+1}\Gamma(\frac{n}{2})
+O\left(\eta^{\eta+\frac{3}{2}}\right)}{(n+1)\Gamma(\frac{n+1}{2})
\Gamma(n)|\det
T|^2}\notag\\
&=
\dfrac{2^{2n-2}\pi^{n-\frac12}\eta^{n+1}\Gamma(\frac{n}{2})
+O\left(\eta^{\eta+\frac{3}{2}}\right)}{(n+1)\Gamma(\frac{n+1}{2})
\Gamma(n)\left(\frac{\sigma_Z}{s_Z}(0)\right)^{n+1}}.\notag
\end{align*}

Using this estimate as before we find that
\begin{equation*}\label{E:hle}
C_n\limsup_{\delta\searrow0}
\frac{\vol( \Omega_{M,\delta})}{\delta^{1/(n+1)}}\le\int_{\bndry\Omega}
\sigma_{\bndry\Omega}.
\end{equation*}

Quoting as before  the existence of peaking functions $h_p$ based on the
Levi polynomial we get the complementary estimate
\begin{equation*}\label{E:hue}
C_n \liminf_{\delta\searrow0}
\frac{\vol( \Omega_{M,\delta})}{\delta^{1/(n+1)}}\ge\int_{\bndry\Omega}
\sigma_{\bndry\Omega}.
\end{equation*}

Combining the estimates we have \eqref{E:mt}. \qed

\bigskip

\section{Comments}\label{E:com}

\begin{enumerate}

\item  The proof of Theorem \ref{T:mt} can easily be adapted to yield the following
result:

\begin{Thm}\label{T:mt2}
Let $\Omega\subset \C^n$ be a relatively compact strictly pseudoconvex domain with
$C^{3}$  boundary inside a complex manifold equipped with a smooth positive 2-form
$\omega$.  Let $\Omega_{M,\delta,\omega}$ be defined as in \eqref{E:cfldef}.  Then
\begin{equation}\label{E:mt2}
C_n
\lim_{\delta\searrow0}
\frac{\vol_\omega( \Omega_{M,\delta,\omega})}{\delta^{1/(n+1)}}=\int_{\bndry\Omega}
\sigma_{\bndry\Omega,\omega}.
\end{equation}
\end{Thm}

\item It would be interesting to know if a result like Theorem \ref{T:mt} holds
also for weakly pseudoconvex domains (possibly involving some
reformulation of the family $P_M(\Omega)$), and if a limit like
\eqref{E:mt} can be shown to exist (independent of the choice of $M\ge
M_0$) also in non-smooth settings.

Note that in the case of a polydisk, the limit \eqref{E:mt} does exist and in fact
it vanishes. 

\medskip

\item For a general bounded convex body $K$ in $\R^n$ it is known that the
limit
\eqref{E:swlim}  coincides with the integral 
\begin{equation}\label{E:eas}
\int_{\bndry K}\kappa_{\bndry K}^{1/(n+1)}  s_{\bndry K},
\end{equation}
where  $s_{\bndry K}$ denotes $(n-1)$-dimensional Hausdorff measure and
$\kappa_{\bndry K}$ denotes the Gaussian curvature of $\bndry\kappa$
 which exists in a suitable pointwise sense almost everywhere with respect to
$s_{\bndry K}$ (so in essence the singular part of the curvature is
discarded). (See [ScWe], [Lei1,
\S2.7].)

In the holomorphic setting there is no evident way to similarly interpret
Fefferman's measure on boundaries of arbitrary bounded pseudoconvex
domains.  But if we impose additional hypotheses such an interpretation may
be possible.  It would be interesting to know if this can be carried out in
particular for domains with Lipschitz boundary which admit strictly
plurisubharmonic defining functions.  

\medskip

\item A number of results have been proved relating affine surface area to
the complexity of approximating polytopes (see the survey [Gru]).  It would
be interesting to have similar results in the holomorphic setting concerning 
approximation by analytic polyhedra.  (Some natural-sounding
notions of complexity of analytic polyhedra definitely will not work for this
purpose: see [Hrm, Lemma 5.3.8].)
 
\end{enumerate}


\end{document}